\pdfoutput=1
\documentclass[twoside]{article}
\makeatletter
\def\ps@myheadings{%
     \let\@oddfoot\@empty\let\@evenfoot\@empty
     \def\@evenhead{\thepage\hfil\MakeUppercase{\footnotesize\leftmark}\hfil}%
     \def\@oddhead{{\hfil\MakeUppercase{\footnotesize\rightmark}}\hfil\thepage}%
     \let\@mkboth\@gobbletwo
     \let\sectionmark\@gobble
     \let\subsectionmark\@gobble
     }
\makeatother
\usepackage{amsmath}
\usepackage{amsthm}
\usepackage{amscd}
\usepackage{amssymb}
\usepackage{amsfonts}
\usepackage{amscd}
\usepackage{amsbsy}
\usepackage{mathrsfs}
\usepackage{comment}

\usepackage{epsf}
\usepackage{epsfig}
\usepackage{psfig}
\usepackage{graphicx}

\newtheorem{theorem}{Theorem}[section]
\newtheorem{lemma}[theorem]{Lemma}
\newtheorem{prop}[theorem]{Proposition}

\theoremstyle{definition}

\theoremstyle{remark}
\newtheorem{remark}[theorem]{Remark}

\numberwithin{equation}{section}

\newcommand{\al}{\alpha}

\newcommand{\Ga}{\Gamma}
\newcommand{\ga}{\gamma}
\newcommand{\sa}{\sigma}
\newcommand{\Sa}{\Sigma}

\newcommand{\La}{\Lambda}

\newcommand{\Om}{\Omega}

\newcommand{\scr}{{\mathcal R}}

\newcommand{\p}{\partial}

\newcommand{\ra}{\rightarrow}

\newcommand{\bs}{\backslash}

\newcommand{\col}{\!:\!}

\newcommand{\fix}{\operatorname{fix}}

\newcommand{\card}{\operatorname{card}}

\newcommand{\is}{\operatorname{Isom}}

\newcommand{\Mob}{\mathop{\rm M\ddot{o}b}\nolimits}

\newenvironment{pf}{\begin{trivlist}\item[]{\bf Proof:\ }}
{\mbox{}\hfill\rule{.08in}{.08in}\end{trivlist}}

\newenvironment{cnstr}{\begin{trivlist}\item[]{\bf Construction:\ }}
{\mbox{}\hfill\rule{.18in}{.18in}\end{trivlist}}

\begin{document}

\pagestyle{myheadings}
\markboth{Boris Apanasov}{Bounded Locally Homeomorphic Quasiregular Mappings}
\title{Hyperbolic topology and bounded locally homeomorphic quasiregular mappings in 3-space}
\author{Boris N. Apanasov}

\date{}

\maketitle

\begin{abstract}
We use our new type of bounded locally homeomorphic quasiregular mappings in the unit 3-ball to address long standing problems for such mappings. The construction of such mappings comes from our construction of non-trivial compact 4-dimensional cobordisms $M$ with symmetric boundary components and whose interiors have complete 4-dimensional real hyperbolic structures. Such bounded locally homeomorphic quasiregular mappings are defined in the unit 3-ball $B^3\subset \mathbb{R}^3$ as mappings equivariant with the standard conformal action of uniform hyperbolic lattices $\Gamma\subset \operatorname{Isom} H^3$ in the unit 3-ball and with its discrete representation $G=\rho(\Gamma)\subset \operatorname{Isom} H^4 $. Here $G$ is the fundamental group of our non-trivial hyperbolic 4-cobordism $M=(H^4\cup\Omega(G))/G$  and the kernel of the homomorphism $\rho\!:\! \Gamma\rightarrow G$ is a free group $F_3$ on three generators.
\end{abstract}

\maketitle
    \footnotetext[2]{
    \hfil\hfil\hfil}

\section{Introduction}
The theory of quasiregular mappings, initiated by the works of M.A.Lavrentiev and later by Reshetnyak
and Martio, Rickman and V\"{a}is\"{a}l\"{a}, shows that they form (from the geometric
function theoretic point of view) the correct generalization of the class of analytic
functions to higher dimensions. In particular, Reshetnyak proved that non-constant
quasiregular mappings are (generalized) branched covers, that is, continuous, discrete and open mappings and hence local homeomorphisms modulo an exceptional set of (topological) codimension at least two, and that they preserve sets of measure zero. For the theory of quasiregular mappings, see the monographs \cite{Re}, \cite{Ri} and \cite{Vu3}.

Here we address some properties of bounded quasiregular mappings $f\,:\, B^3\rightarrow \mathbb{R}^3$ in the unit ball $B^3$ and well known problems on such quasiregular mappings. These results will be heavily based on our recent construction Apanasov \cite{A7} (Theorem 4.1) of surjective locally homemorphic quasiregular mappings
$F\,:\, S^3\backslash S_*\rightarrow S^3$ with topological barriers at points of a dense subset $S_*\subset S^2\subset S^3$. Due to its importance for understanding of results of this paper, in the Appendix we will give some details of this construction based on properties of non-trivial compact 4-dimensional cobordisms $M^4$ with symmetric boundary components - see \cite{A5}-\cite{A7}. The interiors of these 4-cobordisms have complete 4-dimensional real hyperbolic structures and universally covered by the real hyperbolic space $H^4$, while the boundary components of $M^4$ have (symmetric) 3-dimensional conformally flat structures obtained by deformations of the same hyperbolic 3-manifold whose fundamental group $\Gamma$ is a uniform lattice in $\is H^3$. Such conformal deformations of hyperbolic manifolds are well understood after their discovery in \cite{A1}, see \cite{A4}. Nevertheless till recently such "symmetric" hyperbolic 4-cobordisms were unknown despite our well known constructions of non-trivial hyperbolic homology 4-cobordisms with very assymmetric boundary components - see \cite{AT} and \cite{A2}-\cite{A4}.

The above subset $S_*$ of the boundary sphere $S^2=\partial B^3$ is a countable $\Gamma$-orbit of a Cantor subset with Hausdorff dimension $\ln 5/\ln 6 \approx 0.89822444$ (where a uniform hyperbolic lattice $\Gamma$ conformally acts in the unit ball $B^3$). All its points are essential singularities of the bounded locally homeomorphic quasiregular mapping $f\,:\, B^3\rightarrow \mathbb{R}^3$ defined as the restriction to the unit ball $B^3$ of our quasiregular mapping $F\,:\, S^3\backslash S_*\rightarrow S^3$. This bounded quasiregular mapping $f\,:\, B^3\rightarrow \mathbb{R}^3$ has no radial limits at all points $x\in S_*\subset S^2=\partial B^3$ and gives an advance to the well known Pierre Fatou's problem on the correct analogue for higher-dimensional quasiregular mappings of the Fatou's theorem \cite{F} on radial limits of a bounded analytic function of the unit disc. There are several results concerning radial limits of mappings of the unit ball. The most recent progress is due to Kai Rajala who in particular proved that radial limits exist for infinitely many points of the unit sphere, see \cite{Ra} and references there for some earlier results in this direction.

Another application of our construction Apanasov \cite{A7} (Theorem 4.1) of surjective locally homemorphic quasiregular mappings $F\,:\, S^3\backslash S_*\rightarrow S^3$
is to a well known open problem on injectivity of quasiregular mappings in space formulated by Matti Vourinen in 1970-1980s (see Vuorinen \cite{Vu1}-\cite{Vu2},\cite{Vu3} (page 193, Problem 4) and Problem 7.66 in the Hayman's list of problems \cite{BBH}, \cite{HL}). The problem asks whether a proper quasiregular mapping $f$ in the unit ball $B^n$, $n\geq 3$, with a compact branching set $B_f\subset B^n$ is injective. It is false when $n=2$. The conjecture is known to be true in the special case $f(B^n)=B^n$, $n\geq 3$ - see Vuorinen \cite{Vu1}. In Section 2 we show that the quasiregular mapping $f\,:\, B^3\rightarrow \mathbb{R}^3$ defined as the restriction to the unit ball $B^3$ of our quasiregular mapping $F\,:\, S^3\backslash S_*\rightarrow S^3$ from Apanasov \cite{A7} (Theorem 4.1) is essentially a counter-example to this conjecture. This mapping $f$ is essentially proper in the sense that any compact $C\subset f(B^3)$ has a compact subset $C'\subset B^3$ such that $f(C')=C$. This mapping $f$ is bounded locally homeomorphic but not injective quasiregular mapping in the unit ball. Its restriction to a round ball $B_r\subset B^3$ of radius $r<1$ arbitrary close to one gives (after re-scaling) a proper bounded quasiregular mapping of the round ball $B^3$ serving as a counter-example to this Vuorinen conjecture (Theorem \ref{inj}).

The last task of this paper is to investigate the asymptotic behavior of bounded locally homeomorphic quasiregular mappings in the unit ball in smaller balls $B^n(r)\subset B^n$ of radius $r$ close to one. There is an open Matti Vuorinen conjecture that in dimension $n\geq 3$ it is not possible that for $y \in f(B^n)$ and all $r\in (1/2,1)$, the cardinality
\begin{eqnarray}\label{prob3}
 \card (B^n(r) \cap(f^{-1}(y))) > \frac{1}{(1-r)^{n-1}}
\end{eqnarray}

 In Section 3 we show that this question is closely related to the growth function of the kernel (a free group of rank 3) $F_3\subset\Gamma$ of the homomorphism $\rho\col \Gamma \rightarrow G$ of our hyperbolic lattice $\Gamma\subset \is H^3$ to the constructed discrete group $G\subset \is H^4$ - see \cite{A7}. 
 We conclude that our bounded locally homeomorphic quasiregular mappings in the unit ball $B^3$ satisfies this conjecture.

 \subsection{Acknowledgments}

The author is grateful to Matti Vuorinen for attracting our attention to these problems and for fruitful discussions.

 \section{Not injective bounded quasiregular mappings}

 Here we apply our construction \cite{A7} (see Appendix: Theorems \ref{constr} and \ref{map}) of bounded locally homeomorphic quasiregular mapping $F\col B^3\rightarrow \mathbb{R}^3$ to solve the Matti Vourinen open problem on injectivity of quasiregular mappings in 3-dimensional space. This well known problem was formulated by Matti Vourinen in 1970-1980s as result of investigations of quasiregular mappings in space (see Vuorinen \cite{Vu1}-\cite{Vu2},\cite{Vu3} (page 193, Problem 4) and Problem 7.66 in the Hayman's list of problems \cite{BBH}, \cite{HL}).

 The problem asks whether a proper quasiregular mapping $f$ in the unit ball $B^n$, $n\geq 3$, with a compact branching set $B_f\subset B^n$ is injective. The mapping $f(z)=z^2$, where $z\in B^2$, shows that the conjecture is false when $n=2$. The conjecture is known to be true in the special case $f(B^n)=B^n$, $n\geq 3$ - see Vuorinen \cite{Vu1}. Here we give a counter-example to this conjecture for $n=3$.

 \begin{prop}\label{e-inj}
 Let the uniform hyperbolic lattice $\Gamma\subset\is H^3$ and its discrete representation $\rho\col\Gamma\rightarrow G\subset\is H^4$ with the kernel as a free subgroup $F_3\subset\Gamma$ be as in Theorem \ref{constr}. Then the bounded locally homeomorphic quasiregular mapping $F\col B^3\rightarrow \mathbb{R}^3$ constructed in Theorem \ref{map} as a $\Gamma$-equivariant mapping in (\ref{QR}) is an essentially proper bounded quasiregular mapping in the unit 3-ball $B^3$ which is locally homeomorphic ($B_F=\emptyset$) but not injective.
\end{prop}

 \begin{pf} The discrete group $G=\rho(\Gamma)\subset\is H^4\cong\Mob (3)$ has its invariant bounded connected component $\Omega_1\subset\Omega(G)\subset S^3$ where its fundamental polyhedron $P_1$ is quasiconformally homeomorphic to the convex hyperbolic polyhedron $P_0$ fundamental for our hyperbolic lattice $\Gamma\subset\is H^3$ conformally acting in the unit ball $B^3(0,1)$), $\phi_1^{-1}\col P_0\rightarrow P_1$. This homeomorphism $\phi_1^{-1}$ maps polyhedral sides of $P_0$ to the corresponding sides of the polyhedron $P_1$ and preserves all dihedral angles.

 Our bounded locally homeomorphic quasiregular mapping $F\col B^3\rightarrow \Omega_1\subset \mathbb{R}^3$ defined in (\ref{QR}) as the equivariant extention of this homeomorphism $\phi_1^{-1}$ maps the tessellation of $B^3$ by compact $\Gamma$-images of $P_0$ to the tessellation of $\Omega_1$ by compact $G$-images of $P_1$. This shows that for any compact subset $C\subset\Omega_1=F(B^3)$ (covered by finitely many polyhedra $g(P_1)$, $g\in G$), there is a compact subset $C'\subset B^3$ (covered by finitely many corresponding polyhedra $\ga(P_1)$, $\ga\in \Ga$) such that $F(C')=C$.

 On the other hand, this locally homeomorphic quasiregular mapping $F$ is not injective in $B^3$. In fact, for any element $\gamma\neq 1$ in the kernel $F_3\subset\Gamma$ of the homomorphism $\rho\col\Gamma\rightarrow G$ the image $F(P_0)=P_1$ of the fundamental polyhedron $P_0$ of the lattice $\Gamma$ is the same
as the image $F(\gamma(P_0))$ of the translated polyhedron $\gamma(P_0)$, $P_0\cap \gamma(P_0)=\emptyset$.
 \end{pf}
One may restrict our not injective essentially proper bounded quasiregular mapping $F$ in the unit 3-ball $B^3$ from Proposition \ref{e-inj} to a round ball $B_r\subset B^3$ of radius $r<1$ arbitrary close to one. The composition of this restriction $F_r\col B_r \rightarrow \mathbb{R}^3$ with the stretching of $B_r$ to the unit ball $B^3$, i.e. the mapping

\begin{eqnarray}\label{prob2}
 f\col B^3 \rightarrow \mathbb{R}^3\,, f(x)=F_r(rx)
\end{eqnarray}
is a proper bounded quasiregular mapping of the unit ball $B^3$. For any point $y\in f(B^3)$ the number $N_y$ of its pre-images, $N_y=\mid\{x\in B^3\col f(x)=y\}\mid$, is finite
(here $\mid E \mid$ denotes the cardinality of the set $E$). The number $N_y$ of such pre-images of $y\in f(B^3)$ is determined by the number of images $\gamma(P_{ker})$, $\ga\in F_3\subset\Gamma$,
of a fundamental polyhedron $P_{ker}\subset H^3$ in our round ball $B_r\subset B^3$ of radius $r<1$ defining the mapping $f$ in (\ref{prob2}). Here $F_3$ is the free subgroup $F_3\subset\Gamma$ in the uniform hyperbolic lattice $\Ga\subset \is H^3$ (the kernel of the discrete representation $\rho\col\Gamma\rightarrow G\subset\is H^4$ from Theorems \ref{e-inj} and \ref{constr}), and $P_{ker}\subset H^3$ is its fundamental polyhedron in the hyperbolic space $H^3$. Making the radius $r<1$ sufficiently close to $1$, one can make the number $N_y$ arbitrary large.  This proves the following (the Vuorinen conjecture' counter-example):

\begin{theorem}\label{inj}
There are proper bounded quasiregular mappings $f\col B^3 \rightarrow \mathbb{R}^3$ without branching sets ($B_f=\emptyset$) which are locally homeomorphic but not injective. Their pre-images $\{x\in B^3\col f(x)=y\}$ are finite and can be made arbitrary large.
\end{theorem}

\section{Asymptotics of bounded quasiregular mappings in the unit ball and growth in free groups}

Here we investigate the asymptotic behavior of bounded locally homeomorphic quasiregular mappings $f$ in the unit ball. The question is how many pre-images of a point $y \in f(B^n)$ do we have in smaller balls $B^n(r)\subset B^n$ of radius $r$ close to one. There is an open Matti Vuorinen conjecture that in dimension $n\geq 3$ it is not possible that for $y \in f(B^n)$ and all $r\in (1/2,1)$, the cardinality of such pre-image in $B^n(r)$ is bigger than $(1-r)^{1-n}$ - see (\ref{prob3}).

As we show this question for our bounded quasiregular mappings in the unit ball $B^3$, $F\col B^3\rightarrow \mathbb{R}^3$, constructed in Theorem \ref{map} is closely related to the growth function of the kernel $F_3\subset\Gamma$ of the homomorphism $\rho\col \Gamma \rightarrow G$ of our uniform hyperbolic lattice $\Gamma\subset \is H^3$ to the constructed discrete group $G\subset \is H^4$ - see Proposition \ref{homo} and Lemma \ref{ker} in Appendix. Here $F_3$ is a free group on 3 generators.

For free groups $F_{m}$  on $m$ generators one can use well known facts about their growth functions, cf. \cite{Gri}. The growth function $\ga_{G, \Sigma}$ of a group $(G, \Sigma)$  with a generating set $\Sigma$ counts the number of elements in $G$ whose length (in the word metric) is at most a natural number $k$:
\begin{equation}
\ga_{G, \Sigma}(k)=\mid \{g\in G\col |g|_{\Sigma}\leq k\}\mid
\end{equation}
where $\mid E \mid$ denotes the cardinality of the set $E$, and $k$ is a natural number.

\begin{lemma}\label{free-growth}
A free group $F_m$ on $m$ generators (for any free system $\Sigma$ of generators) has $2m(2m-1)^{k-1}$ elements of length $k$, and its growth function:
\begin{equation}
\ga_{F_m}(k)=1+\frac{m}{m-1}((2m-1)^k-1).
\end{equation}
\end{lemma}

\begin{pf} Clearly in a free group $F_m$ on $m$ generators we have the number of elements with length $i$ equals to $c_i=\mid \{g\in F_m\col |g|=i\}\mid=2m(2m-1)^i$. Therefore the growth function $\ga_{F_m}(k)=1+2m+2m(2m-1)+ \ldots +2m(2m-1)^{k-1}$. This gives the value $\ga_{F_m}(k)$ in the Lemma.
\end{pf}

 In the embedded Cayley graph
$\varphi(K(\Ga,\Sa))\subset B^3$ (i.e. the graph that is dual to the tessellation of $B^3$ by convex hyperbolic
polyhedra $\ga(P_0)$, $\ga\in \Ga$), we consider its subgraph (a tree) corresponding to our free group $F_3\subset\Ga$ on 3 generators (the kernel of the homomorphism $\rho$).
The embedding $\varphi$ of the Cayley graph $K(\Ga,\Sa)$ is a $\Ga$-equivariant proper embedding. It is a pseudo-isometry, i.e.
for the word metric $(\ast,\ast)$ on $K(\Ga,\Sa)$ and the hyperbolic metric $d$ in the unit ball $B^3$,
there are positive constants $K$ and $K'$ such that $(a,b)/K\leq d(\varphi(a),\varphi(b))\leq K\cdot(a,b)$ for all $a,b\in K(\Ga,\Sa)$ satisfying one of the following two conditions: either $(a,b)\geq K'$ or $d(\varphi(a),\varphi(b))\geq K'$.

Let $D$ be the maximum of hyperbolic length of generators of the kernel $F_3\subset \Ga$.
 All vertices of our tree subgraph corresponding to elements in $F_3$ of length at most $k$ are in the hyperbolic ball centered at $0\in B^3$ with radius $R=Dk$.  This hyperbolic ball corresponds to the Euclidean ball $B^3(0, r)\subset B^3(0,1)$ of radius $r=(e^R-1)/(e^R+1)$.

Multiplying (\ref{prob3}) by $(1-r)^{n-1}$, we see that we need to estimate the asymptotics of
\begin{equation}\label{asymp}
(1-r)^{n-1} \card (B^n(r) \cap(F^{-1}(y))).
\end{equation}
for arbitrary small $\epsilon=(1-r)$, or for arbitrary large $\lambda=\ln ((2/(1-r))-1)$.

In the case of our free group $F_3$ on 3 generators (the kernel of the homomorphism $\rho$), Lemma \ref{free-growth} shows that the growth function $\ga_{F_3}(k)=1+3(5^k-1)/2$. This reduces the asymptotics of (\ref{asymp}) to the asymptotics of $3(5^{\lambda/D}-1)/e^{2\lambda}$ for arbitrary large $\lambda$. Since the last expression tends to $0$ when $\lambda$ tends to $\infty$,
we conclude that our bounded locally homeomorphic quasiregular mappings $F\col B^3\rightarrow \Omega_1\subset\mathbb{R}^3$ in the unit ball $B^3$ satisfy the Vuorinen conjecture (\ref{prob3}).

\begin{remark}
There is an important observation. If in our analysis of the asymptotics of (\ref{asymp}) (and in our construction of groups $\Ga$ and $G$) the kernel of the corresponding homomorphism $\rho\col\Gamma\rightarrow G\subset\is H^4$ were a free subgroup $F_m$ on a big number $m$ of generators, then our last expression would tend to $\infty$ when $\lambda$ tends to $\infty$. This would provide a way to constructing a similar bounded locally homeomorphic quasiregular mapping in the unit ball giving a possible counter example to (\ref{prob3}).
\end{remark}

\section{Appendix: Hyperbolic 4-cobordisms and deformations of hyperbolic structures}

For the readers convenience, here we provide essential details of our construction \cite{A7} of locally homeomorphic quasiregular mappings $F\col S^3\setminus S_*\rightarrow S^3$ based on the properties of non-trivial "symmmetric" hyperbolic 4-cobordisms $M^4=(H^4\cup\Omega(G))/G$ constructed in Apanasov \cite{A5}. Properties of the fundamental group $\pi_1(M^4)\cong G\subset \is H^4$ of such "symmmetric" hyperbolic 4-cobordisms $M^4=H^4/G$  acting discretely in the hyperbolic 4-space $H^4$ and in the discontinuity set $\Omega(G)\subset \partial H^4 = S^3$ are very essential for our construction of the quasiregular mapping $F$.

We start with our construction \cite{A5} of such discrete group $G\subset \is H^4$ and the corresponding discrete representation $\rho\col \Gamma \rightarrow G$ of a uniform hyperbolic lattice $\Gamma\subset \is H^3$.
These discrete groups $G$ and $\Gamma$ produce a non-trivial (not a product) hyperbolic 4-cobordisms $M^4=(H^4\cup\Omega(G))/G$ whose boundary components $N_1$ and $N_2$ are topologically and geometrically symmetric to each other. These $N_1$ and $N_2$ are covered by two $G$-invariant connected components $\Omega_1$ and $\Omega_2$ of the discontinuity set $\Omega (G)\subset S^3$, $\Omega(G)=\Omega_1\cup\Omega_2$. The conformal action of $G=\rho(\Gamma)$ in these components $\Omega_1$ and $\Omega_2$ is symmetric and has contractible fundamental polyhedra $P_1\subset \Omega_1$ and $P_2\subset \Omega_2$ of the same combinatorial type allowing to realize them as a compact polyhedron $P_0$ in the hyperbolic 3-space, i.e. the dihedral angle data of these polyhedra satisfy the Andreev's conditions \cite{An1}. Nevertheless this geometric symmetry of boundary components of our hyperbolic 4-cobordism $M^4(G)$) does not make the group $G=\pi_1(M^4)$ quasi-Fuchsian, and our 4-cobordism $M^4$ is non-trivial.

 Here a Fuchsian
group $\Ga\subset\is H^3\subset \is H^4$ conformally acts in the 3-sphere $S^3=\p H^4$ and preserves a round ball $B^3\subset S^3$ where it acts as a cocompact discrete group of isometries of $H^3$. Due to the Sullivan
structural stability (see Sullivan \cite{S} for $n=2$ and Apanasov \cite{A2}, Theorem 7.2), the
space of quasi-Fuchsian representations of a hyperbolic lattice $\Ga\subset\is H^3$ into $\is H^4$
is an open connected component of the Teichm\"uller space of $H^3/\Ga$ or the variety of conjugacy classes of discrete representations
$\rho\col \Ga\ra\is H^4$. Points in this (quasi-Fuchsian) component of the variety correspond to trivial hyperbolic 4-cobordisms
$M(G)$ where the discontinuity set $\Om(G)=\Om_1\cup\Om_2\subset S^3=\p H^4$ is the union of two topological 3-balls $\Om_i$, $i=1,2$, and
$M(G)$ is homeomorphic to the product of $N_1$ and the closed interval $[0,1]$.

 We may consider hyperbolic 4-cobordisms $M(\rho(\Ga))$ corresponding to uniform hyperbolic lattices $\Ga\subset \is H^3$ generated by reflections. Natural inclusions of these lattices into $\is H^4$ act at infinity $\p H^4 = S^3$ as Fuchsian groups
$\Ga\subset \Mob (3)$ preserving a round ball $B^3\subset S^3$. In this case
 our construction of the mentioned discrete groups $\Gamma$ and $G=\rho(\Gamma)$ results in the following (see Apanasov \cite{A5}):

\begin{theorem}\label{constr}
There exists a discrete M\"obius group $G\subset \Mob (3)$ on the 3-sphere $S^3$ generated by finitely many reflections such that:
\begin{enumerate}
\item Its discontinuity set $\Om(G)$ is the union
of two invariant components $\Om_1$, $\Om_2$;
\item Its fundamental polyhedron $P\subset S^3$ has two contractible components $P_i\subset\Om_i$, $i=1,2$,
having the same combinatorial type (of a compact hyperbolic polyhedron $P_0\subset H^3$);
\item For the uniform hyperbolic lattice $\Ga\subset\is H^3$ generated by reflections in sides of the hyperbolic
polyhedron $P_0\subset H^3$ and acting on the sphere $S^3=\p H^4$ as a discrete Fuchsian group  $i(\Ga)\subset \is H^4=\Mob(3)$ preserving a round ball $B^3$ (where $i\col\is H^3\subset\is H^4$ is the natural inclusion),
the group $G$ is its image under a homomorphism $\rho\col\Ga\ra G$  but it is not quasiconformally (topologically) conjugate in $S^3$ to $i(\Ga)$.
\end{enumerate}
\end{theorem}

\begin{cnstr}
We define a finite collection $\Sa$ of reflecting 2-spheres
$S_i\subset S^3$, $1\leq i \leq N$.
As the first three spheres $S_1, S_2$ and  $S_3$ we take the coordinate planes $\{x\in {\mathbb R}^3\col x_i=0\}$, and
 $S_4=S^2(0,R)$ is the round sphere of some radius $R>0$ centered at the origin. The value of the radius $R$ will be determined later.
Let $B=\bigcup_{1\leq i\leq 4} B_i$ be the union of the closed balls bounded by these four spheres,
and let $\p B$ be its boundary (a topological 2-sphere) having four vertices which are the intersection points of four triples of our spheres.
We consider a simple closed loop $\alpha\subset \p B$
which does not contain any of our vertices and which symmetrically separates two pairs of these vertices from each other
as the white loop does on the tennis ball. This loop
$\alpha$ can be considered as the boundary of a topological 2-disc $\sigma$ embedded in the complement
   $D=S^3\setminus B$ of our four balls. Our geometric construction needs a detailed
   description of such a 2-disc $\sa$ and its boundary loop $\al=\p \sa$ obtained as it is shown in
   Figure \ref{fig7}.

The desired disc $\sa\subset D=S^3\setminus B$ can be described as the boundary in the domain $D$ of the union
of a finite chain of adjacent blocks $Q_i$ (regular cubes) with disjoint interiors whose centers
lie on the coordinate planes $S_1$ and $S_2$ and whose sides are parallel to the coordinate planes.
This chain starts from the unit cube whose center lies
in the second coordinate axis, in $e_2\cdot \mathbb R_{+}\subset S_1\cap S_3$. Then our chain goes up through small adjacent cubes centered in the coordinate plane $S_1$, at some point changes its direction to the horizontal one toward the third coordinate axis, where it turns its horizontal direction by a right angle again (along the coordinate plane $S_2$), goes toward the vertical line passing through the second unit cube centered in
$e_1\cdot \mathbb R_{+}\subset S_2\cap S_3$, then goes down along that vertical line and finally ends
at that second unit cube, see Figure \ref{fig7}. We will define the size of small cubes $Q_i$ in our block chain
and the distance of the centers of two unit cubes to the origin in the next step of our construction.

\begin{figure}
\centering
\epsfxsize=10cm
\epsfbox{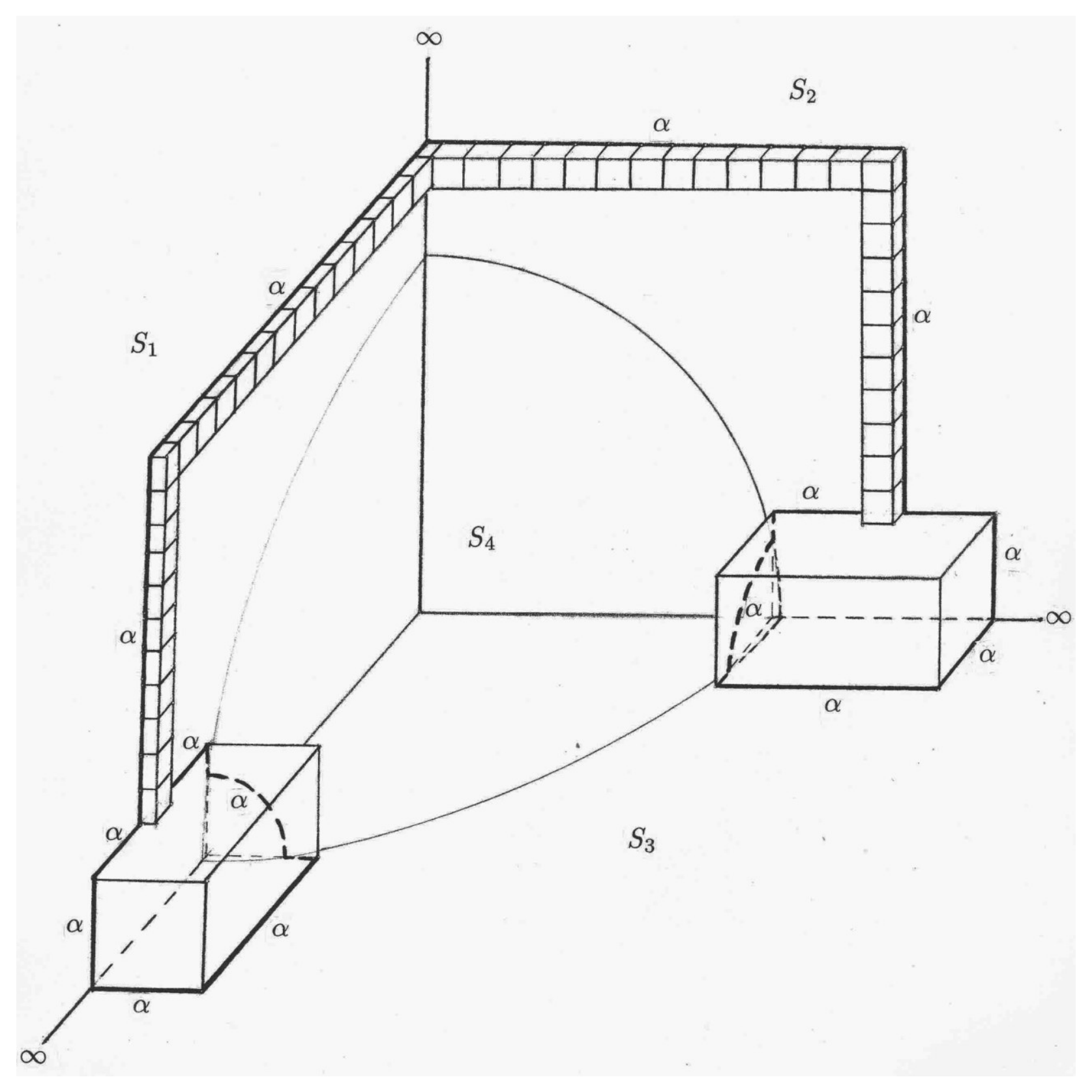}
\caption{Configuration of blocks and the white loop $\alpha\subset \p B$.}
\label{fig7}
\end{figure}

\begin{figure}
\centering
\epsfxsize=7cm
\epsfbox{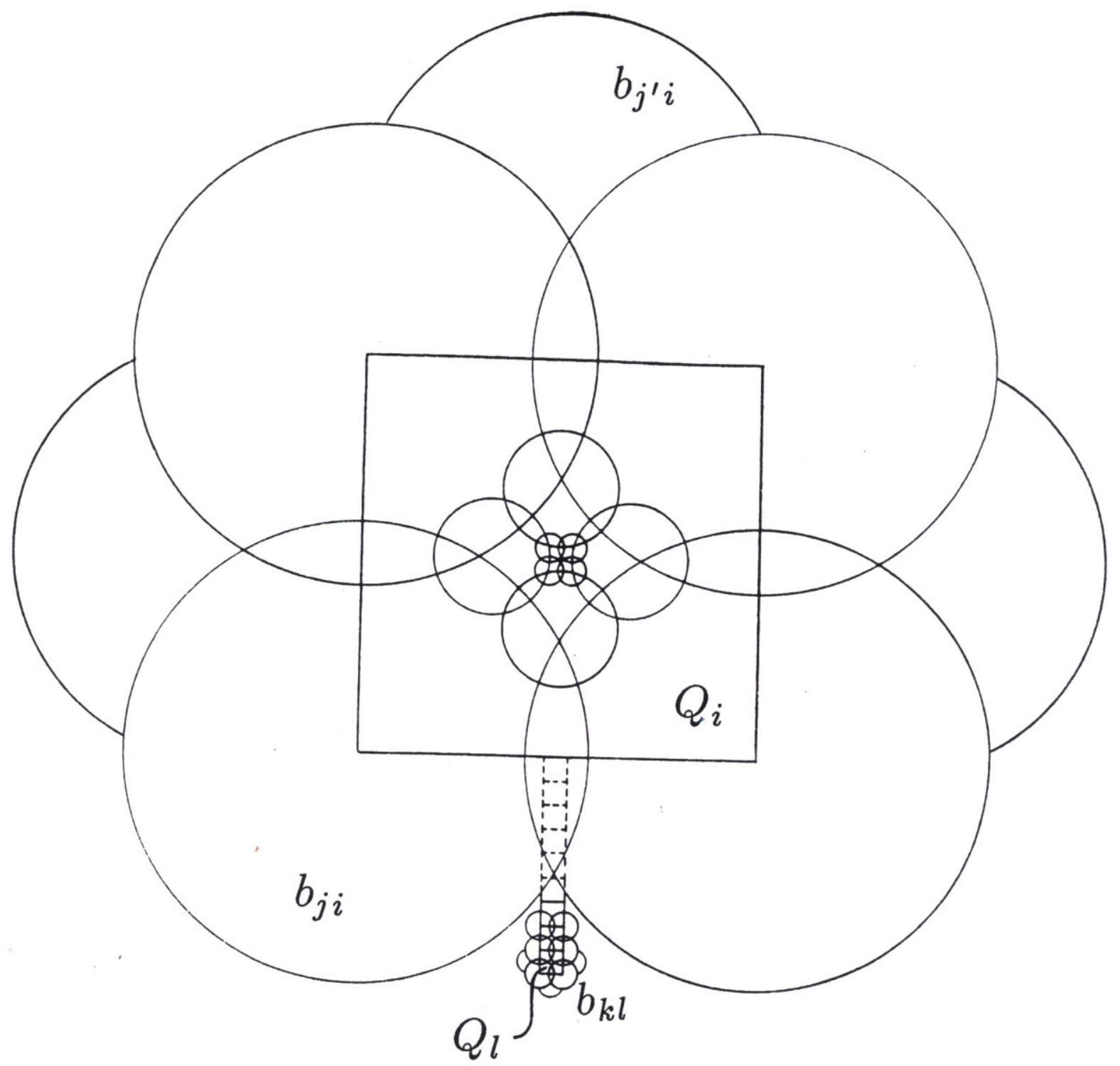}
\caption{Big and small cube sizes and ball covering}
\label{fig4}
\end{figure}

Let us consider one of our cubes $Q_i$, i.e. a block of our chain, and let $f$ be its square side having a nontrivial intersection with our 2-disc
$\sa\subset D$.
For that side $f$ we consider spheres $S_j$ centered at its vertices and
having a radius such that each two spheres centered at the ends of an edge of $f$ intersect each other with angle $\pi/3$.
In particular, for the unit cubes such spheres have radius $ \sqrt{3}/3 $. From such defined spheres we select those spheres that have centers
in our domain $D$ and then include them in the collection $\Sa$ of reflecting spheres.
Now we define the distance of the centers of our big (unit) cubes to the origin. It is determined by the condition that the sphere $S_4=S^2(0,R)$
is orthogonal to the sphere $S_j\in\Sa$ centered at the vertex of such a cube closest to the origin.
As in Figure \ref{fig4}, let $f$ be a square side of one of our cubic blocks  $Q_i$ having a nontrivial intersection $f_{\sa}=f\cap\sa$ with our
2-disc $\sa\subset D$.
We consider a ring of four spheres $S_i$ whose centers are interior points of $f$ which lie outside of the four
previously defined spheres $S_j$ centered at vertices of $f$ and such that each sphere $S_i$ intersects two adjacent spheres $S_{i-1}$ and $S_{i+1}$
(we numerate spheres $S_i$ mod 4) with angle $\pi/3$. In addition these spheres $S_i$
are orthogonal to the  previously defined ring of bigger spheres $S_j$, see Figure \ref{fig4}.
From such defined spheres $S_i$ we select those spheres that have nontrivial
intersections with our domain $D$ outside the previously defined spheres $S_j$, and then include them in the collection $\Sa$ of reflecting spheres.
If our side $f$ is not the top side of one of the two unit cubes we add another sphere $S_k\in\Sa$.
It is centered at the center of this side $f$ and is orthogonal to the four previously defined spheres $S_i$ with centers in $f$, see Figure \ref{fig4}.

Now let $f$ be the top side of one of the two unit cubes of our chain. Then, as before,
we consider another ring of four spheres $S_k$. Their centers are interior points of $f$, lie outside of the four
previously defined spheres $S_i$ closer to the center of $f$ and such that each sphere $S_k$ intersects two adjacent spheres $S_{k-1}$ and $S_{k+1}$
(we numerate spheres $S_k$ mod 4) with angle $\pi/3$. In addition these new four spheres $S_k$
are orthogonal to the previously defined ring of bigger spheres $S_i$, see Figure \ref{fig4}.
 We note that the centers of these four new spheres $S_k$ are vertices of a small square
$f_s\subset f$ whose edges are parallel to the edges of $f$, see Figure \ref{fig4}. We set this square $f_s$ as the bottom side of the small cubic box adjacent to the unit one.
This finishes our definition of the family of twelve round spheres whose interiors cover the square ring $f\bs f_s$ on the top side of one of the two unit cubes
in our cube chain and tells us which two spheres among the four new defined spheres $S_k$ were already included
in the collection $\Sa$ of reflecting spheres
(as the spheres $S_j\in\Sa$ associated to small cubes in the first step).

This also defines the size of small cubes in our block chain. Now we can vary the remaining free parameter $R$
(which is the radius of the sphere $S_4\in\Sa$) in order to make two horizontal rows of small blocks with centers in $S_1$ and $S_2$, correspondingly,
to share a common cubic block centered at a point in $e_3\cdot \mathbb R_{+}\subset S_1\cap S_2$, see Figure \ref{fig7}.

We can use the collection $\Sa$ of reflecting spheres $S_i$ to define a discrete reflection group $G=G_{\Sa}\subset \Mob(3)$. Important properties of $\Sa$ are: (1) the closure of the disc $\sa\subset D$ is covered by balls $B_j$; (2) any two spheres $S_j, S_{j'}\in\Sa$ either are disjoint or intersect with angle $\pi/2$ or  $\pi/3$;
(3) the complement of all balls $B_j$, $1\leq j \leq N$, is the union of two disjoint contractible
polyhedra $P_1$ and $P_2$ of the same combinatorial type and equal corresponding dihedral angles.
So the discontinuity set $\Omega(G)\subset S^3$ of $G$ consists of two invariant connected components $\Om_1$ and $\Om_2$ which are the unions of the $G$-orbits of $\bar{P_1}$ and $\bar{P_2}$, and this defines a Heegaard splitting of the 3-sphere $S^3$ (see \cite{A7}):
\begin{lemma}\label{h-body}
The splitting of the discontinuity set $\Om\subset S^3$ of our discrete reflection group $G=G_{\Sa}\subset \Mob(3)$ into $G$-invariant components $\Om_1$ and $\Om_2$ defines a Heegaard splitting of the 3-sphere $S^3$ of infinite genus with ergodic word hyperbolic group $G$ action on the separating boundary $\La(G)$.
\end{lemma}

To finish our construction in Theorem \ref{constr} we notice that the combinatorial type (with magnitudes of dihedral angles) of the bounded component $P_1$ of the fundamental polyhedron $P\subset S^3$ coincides with
the combinatorial type of its unbounded component $P_2$. Applying Andreev's theorem on 3-dimensional hyperbolic polyhedra \cite{An1},
one can see that there exists a compact hyperbolic polyhedron $P_0\subset H^3$ of the same combinatorial type with the same dihedral angles ($\pi/2$ or $\pi/3$).
So one can consider a uniform hyperbolic lattice $\Ga\subset\is H^3$ generated by reflections in sides of the hyperbolic
polyhedron $P_0$. This hyperbolic lattice $\Ga$ acts in the sphere $S^3$ as a discrete co-compact Fuchsian group $i(\Ga)\subset \is H^4=\Mob(3)$
(i.e. as the group
$i(\Ga)\subset \is H^4$ where $i\col\is H^3\subset\is H^4$ is the natural inclusion)
preserving a round ball $B^3$ and having its boundary sphere $S^2=\p B^3$ as the limit set. Obviously there is no self-homeomorphism of the sphere $S^3$
conjugating the action of the groups $G$ and $i(\Ga)$ because the limit set $\La(G)$ is not a topological 2-sphere. So the constructed group $G$ is not a
quasi-Fuchsian group.
\end{cnstr}

One can construct a natural homomorphism $\rho\col \Ga\ra G$, $\rho\in\scr_3(\Ga)$, between these two Gromov hyperbolic groups $\Ga\subset \is H^3$ and $G\subset \is H^4$
defined by the correspondence between sides of
the hyperbolic polyhedron $P_0\subset H^3$ and reflecting spheres $S_i$ in the collection $\Sa$ bounding the fundamental polyhedra $P_1$ and $P_2$. Then we have:

\begin{prop}\label{homo}
The homomorphism $\rho\in\scr_3(\Ga)$, $\rho\col \Ga\ra G$, in Theorem \ref{constr} is not an isomorphism. Its kernel $\ker(\rho)=\rho^{-1}(e_G)$ is a
free rank 3 subgroup $F_3\lhd\Ga$.
\end{prop}
Its proof (see \cite{A7}, Prop.2.4) is based on the following statement (see  \cite{A7}, Lemma 2.5) in combinatorial group theory:
 \begin{lemma}\label{ker}
Let $A=\langle a_1, a_2 \mid a_1^2,\, a_2^2,\, (a_1a_2)^2\rangle $ $\cong$ $B=\langle b_1, b_2 \mid b_1^2,\, b_2^2,\, (b_1b_2)^2\rangle $ $\cong$
$C=\langle c_1, c_2 \mid c_1^2,\, c_2^2,\, (c_1c_2)^2\rangle \cong \mathbb{Z}_2 \times \mathbb{Z}_2$,
and let $\varphi\col A\ast B\ra C$ be a homomorphism of the free product $A\ast B$ into $C$ such that $\varphi(a_1)=\varphi(b_1)=c_1$ and
 $\varphi(a_2)=\varphi(b_2)=c_2$. Then the kernel $\ker(\varphi)=\varphi^{-1}(e_C)$ of $\varphi$ is a free rank 3 subgroup $F_3\lhd A\ast B$  generated by elements
 $x=a_1b_1$, $y=a_2b_2$ and $z=a_1a_2b_2a_1=a_1ya_1$.
\end{lemma}

\subsection{Bending homeomorphisms between polyhedra}

Here we sketch our construction of a quasiconformal homeomorphism $\phi_1\col P_1\rightarrow P_0$
between the fundamental polyhedron $P_1\subset\Omega_1\subset\Omega(G)\subset S^3$ for the group $G$ action in $\Omega_1$ and the
fundamental polyhedron $P_0\subset B^3$ for conformal action of our hyperbolic lattice $\Gamma\subset\is H^3$ from Theorem \ref{constr}. This mapping $\phi_1$ is a composition of finitely many elementary "bending homeomorphisms". It maps faces to faces, and preserve the combinatorial structure of the polyhedra and their corresponding dihedral angles.

First we observe that to each cube $Q_j, 1\leq j\leq m$, used in the above construction of the group $G$ (see Figure \ref{fig7}), we may associate a round
ball $B_j$ centered at the center of the cube $Q_j$ and such that its boundary sphere is orthogonal to the reflection spheres $S_i$ from our generating family $\Sigma$ whose centers are at vertices of the cube $Q_j$. In particular for the unit cubes $Q_1$ and $Q_m$, the reflection spheres $S_i$ centered at their vertices have radius $\sqrt{3}/3$, so the balls $B_1$ and $B_m$ (whose boundary spheres are orthogonal to those corresponding reflection spheres $S_i$) should have radius $\sqrt{5/12}$. Also we add another extra ball $B^3(0,R)$ (which we consider as two balls $B_0$ and $B_{m+1}$) whose boundary is the reflection sphere
$S^2(0,R)=S_4\in\Sigma$ centered at the origin and orthogonal to the closest reflection spheres $S_i$ centered at vertices of two unit cubes $Q_1$ and $Q_m$. Our different enumeration of this ball will be used when we consider different faces of our fundamental polyhedron $P_1$ lying on that reflection sphere $S_4$.

Now for each cube $Q_j, 1\leq j\leq m$, we may associate a discrete subgroup $G_j\subset G\subset\Mob(3)\cong\is H^4$ generated by reflections in the spheres  $S_i\in\Sigma$ associated to that cube $Q_j$ - see our construction in Theorem \ref{constr}. One may think about such a group $G_j$ as a result of quasiconformal bending deformations (see \cite{A2}, Chapter 5) of a discrete M\"{o}bius group preserving the round ball $B_j$ associated to the cube $Q_j$ (whose center coincides with the center of the cube $Q_j$). As the first step in such deformations, we define two
quasiconformal ``bending'' self-homeomorphisms of $S^3$, $f_1$ and
$f_{m+1}$, preserving the balls $B_1,\ldots,B_m$ and the set
of their reflection spheres $S_i$, $i\neq 4$, and transferring $\p B_0$ and $\p B{m+1}$ into 2-spheres
orthogonally intersecting $\p B_1$ and $\p B_m$ along round circles $b_1$ and $b_{m+1}$,
respectively - see (3.1)  and Figure 6 in \cite{A7}.

 In the next steps in our bending deformations, for two adjacent cubes  $Q_{j-1}$ and  $Q_j$, let  us denote $G_{j-1,j}\subset G$ the subgroup  generated by reflections with respect to the spheres $S_i\subset\Sigma$ centered at common vertices of these cubes. This subgroup preserves the round circle $b_j=b_{j-1,j}=\p B_{j-1}\cap\p B_j$. This shows that our group $G$ is a result of the so called "block-building construction" (see \cite{A2}, Section 5.4) from the block groups $G_j$ by sequential amalgamated products:
\begin{eqnarray}\label{amalgam}
G=G_1\underset{G_{1,2}}*G_2\underset{G_{2,3}}*\cdots\underset{G_{j-2,j-1}}*G_{j-1}\underset{G_{j-1,j}}* G_j\underset{G_{j,j+1}} *\cdots\underset{G_{m-1,m}}*G_m
\end{eqnarray}

The chain of these building balls $\{B_j\}, 1\leq j\leq m$, contains the bounded polyhedron $P_1\subset\Omega_1$.
For each pair $B_{i-1}$ and $B_i$
with the common boundary circle $b_i=\p B_{i-1}\cap\p B_i$, $1\leq i\leq m$, we construct a quasi-conformal
bending homeomorphism
$f_i$ that transfers $B_i\cup B_{i-1}$ onto the ball $B_i$ and which is
conformal in dihedral $\zeta_i$-neighborhoods of the spherical disks
$\p B_i\backslash\overline{B_{i-1}}$ and
$\p B_{i-1}\backslash\overline{B_i}$ - see (3.3) and Figure 7 in \cite{A7}.
In each $i$-th step, $2\leq i\leq m$, we reduce the number of
balls $B_j$ in our chain by one. The composition
$f_{m+1}f_i f_{i-1}\cdots f_2 f_1$ transfers all spheres
from $\Sigma$ to spheres orthogonal to the boundary sphere of our last ball $B_m$ which we renormalize as the unit ball $B(0,1)$. We note that all intersection angles between these spheres do not change.
We define our quasiconformal homeomorphism
$\phi_1\col P_1\rightarrow P_0$ as the restriction of the composition
$f_{m+1}f_mf_{m-1}\cdots f_2f_1$ of our bending
homeomorphisms $f_j$ on the fundamental polyhedron $P_1\subset\Omega_1$.

\subsection{Bounded locally homeomorphic quasiregular mappings}

Now we define bounded quasiregular mappings $F\col B^3\rightarrow \mathbb{R}^3$ as in Theorem 4.1 in \cite{A7}:

\begin{theorem}\label{map}
 Let the uniform hyperbolic lattice $\Gamma\subset\is H^3$ and its discrete representation $\rho\col\Gamma\rightarrow G\subset\is H^4$ with the kernel as a free subgroup $F_3\subset\Gamma$ be as in Theorem \ref{constr}. Then there is a bounded locally homeomorphic quasiregular mapping $F\col B^3\rightarrow \mathbb{R}^3$ whose all singularities lie in an exceptional subset $S_*$ of the unit sphere $S^2\subset \mathbb{R}^3$ and form a dense in $S^2$ $\Ga$-orbit of a Cantor subset with Hausdorff dimension $\ln 5/\ln 6 \approx 0.89822444$. These (essential) singularities create a barrier for $F$ in the sense that at points $x\in S_*$ the map $F$ does not have radial limits.
\end{theorem}

\begin{cnstr}
We construct our quasiregular mapping $F\col B^3\rightarrow \Omega_1=F(B^3)$ in the unit ball $B^3$ by equivariant extention of the above quasiconformal homeomorphism $\phi_1^{-1}\col P_0\rightarrow P_1$ which maps polyhedral sides of $P_0$ to the corresponding sides of the polyhedr $P_1$ and preserves combinatorial structures of polyhedra as well as their dihedral angles:
\begin{eqnarray}\label{QR}
F(x)=\rho(\gamma)\circ \phi_1^{-1}\circ \gamma^{-1}(x)  & \text{if}\ |x|<1,\, x\in\gamma(P_0), \,\gamma\in\Gamma
\end{eqnarray}
The tesselations of $B^3$ and $\Omega_1$ by corresponding $\Gamma$- and $G$-images of their fundamental polyhedra $P_o$ and $P_1$
are perfectly similar. This implies that our quasiregular mapping $F$ defined by (\ref{QR}) is bounded and locally homeomorphic.
It follows from Lemma \ref{h-body} that the limit set $\Lambda(G)\subset S^3$ of the group $G\subset \Mob(3)$ defines a Heegard splitting of infinite genus of the 3-sphere $S^3$ into two connected components $\Om_1$ and $\Om_2$ of the discontinuity set $\Om(G)$. The action of $G$  on the limit set $\La(G)$ is an ergodic word hyperbolic action.
For this ergodic action the set of fixed points of loxodromic elements $g\in G$ (conjugate to similarities in $\mathbb{R}^3$) is dense in $\La(G)$. Preimages $\ga\in\Ga$ of such loxodromic elements $g\in G$ for our homomorphism $\rho\col\Ga\rightarrow G$ are loxodromic elements in $\Ga$ with two fixed points $p,q\in\La(\Ga)=S^2, p\neq q$. This and Tukia's arguments of the group completion (see \cite{T} and \cite{A2}, Section 4.6) show that our mapping $F$ can be continuously extended to the set of fixed points of such elements $\ga\in\Ga$, $F(Fix(\ga))=Fix(\rho(\ga))$. The sense of this continuous extension is that if $\ga\in\Ga$ is a loxodromic preimage of a loxodromic element $g\in G$, $\rho(\ga)=g$, and if $x\in S^3\backslash S^2$ tends to its fixed points $p$ or $q$ along the hyperbolic axis of $\ga$ (in $B(0,1)$ or in its complement $\widehat{B(0,1)}$) (i.e. radially) then $\lim_{|x|\to 1} F(x)$ exists and equals to the corresponding fixed point of the loxodromic element $g=\rho(\ga)\in G$. In that sense one can say that the limit set $\La(G)$ (the common boundary of the connected components $\Om_1, \Om_2 \subset\Om(G)$) is the $F$-image of points in the unit sphere $S^2\subset S^3$. So the mapping $F$ is extended to a map onto the closure $\overline{\Omega_1}=\Omega_1\cup\Lambda(G)\subset\mathbb{R}^3$.

Nevertheless not all loxodromic elements $\ga\in\Ga$ in the hyperbolic lattice $\Ga\subset\is H^3$ have their images $\rho(\ga)\in G$ as loxodromic elements. Proposition \ref{homo} shows that $\ker \rho\cong F_3$ is a free subgroup on three generators in the lattice $\Ga$, and all elements $\ga\in F_3$ are loxodromic. Now we look at radial limits $\lim_{x\to p}F(x)$ when $x$ radially tends to a fixed point $p\in S^2$ of this loxodromic element $\ga\in F_3\subset\Ga$.

Let $K(\Ga,\Sigma)$ be the Cayley graph for a group $\Ga$ with a finite generating set $\Sigma$. Our lattice $\Ga\subset \is H^3$ has an embedding $\varphi$ of its Cayley graph $K(\Ga,\Sigma)$ in the hyperbolic space $H^3\cong B^3$.
For a point $0\in H^3$ not fixed by any $\ga\in\Ga\backslash\{1\}$, vertices $\ga\in K(\Ga,\Sa)$ are mapped to $\ga(0)$, and edges joining vertices
$a,b\in K(\Ga,\Sa)$ are mapped to the hyperbolic geodesic segments $[a(0),b(0)]$.
In other words, $\varphi(K(\Ga,\Sa))$ is dual to the tessellation of $H^3$ by
polyhedra $\ga(P_0)$, $\ga\in \Ga$.  The map $\varphi$ is a $\Ga$-equivariant proper
embedding:  for any compact $C\subset H^3$, its pre-image
$\varphi^{-1}(\varphi(K(\Ga,\Sa))\cap C)$ is compact. Moreover for any convex cocompact group $\Gamma \subset \is H^n$ this embedding $\varphi$ is a pseudo-isometry (see \cite{C} and \cite{A2}, Theorem 4.35), i.e. for the word metric on $K(\Ga,\Sa)$ and the hyperbolic metric $d$,
there are $K>0$ and $K'>0$ such that $|a,b|/K\leq d(\varphi(a),\varphi(b))\leq K\cdot|a,b|$
for all $a,b\in K(\Ga,\Sa)$ such that
either $|a,b|\geq K'$ or $d(\varphi(a),\varphi(b))\geq K'$.

This implies (see \cite{A2}, Theorem 4.38) that the limit set of any convex-cocompact group
$\Ga\subset\Mob (n)$ can be identified with its group completion $\overline{\Ga}$, $\overline{\Ga}=\overline{K(\Ga,\Sa)}\setminus K(\Ga,\Sa)$. Namely there exists a continuous and $\Ga$-equivariant bijection $\varphi_{\Ga}\col \overline{\Ga} \rightarrow \La(\Ga)$.

For the kernel subgroup $F_3=\ker \rho \subset \Ga\subset\is H^3$ and for the above pseudo-isometric embedding $\varphi$, we consider its Cayley subgraph in $\varphi(K(\Ga,\Sa))\subset H^3$ which is a tree - see Figure 5 in \cite {A7}. Since the limit set of $\ker \rho=F_3\subset\Ga$ corresponds to the 'bondary at infinity' $\p_{\infty} F_3$ of $F_3\subset\Ga$ (the group completion $\overline{F_3}$), it is a closed Cantor subset of the unit sphere $S^2$ with Hausdorff dimension $\ln 5/\ln 6 \sim 0.89822444$.

 The $\Ga$-orbit $\Ga(\La(F_3))$ of our Cantor set is a dense subset $S_*$ of $S^2=\La(\Ga)$ because of density in the limit set $\La(\Ga)$ of the $\Ga$-orbit of any limit point. In particular we have such dense $\Ga$-orbit $\Ga(\{p,q\})$ of fixed points $p$ and $q$ of a loxodromic element $\ga\in F_3\subset \Ga$ (the images of $p$ and $q$ are fixed points of $\Ga$-conjugates of such loxodromic elements $\ga\in F_3\subset \Ga$).

On the other hand let $x\in l_{\ga}$ where $l_{\ga}$ is the hyperbolic axis in $B(0,1)$ of an element $\ga\in F_3\subset\Ga$.
Denoting $d_{\ga}$ the translation distance of $\ga$, we have that any segment $[x, \ga(x)]\subset l_{\ga}$ is mapped by our quasiregular mapping $F$ to a non-trivial closed loop $F([x, \ga(x)])\subset\Om_1$, inside of a handle of the handlebody $\Om_1$ (mutually linked with $\Om_2$ - similar to the loops $\beta_1\subset \Om_1$ and $\beta_2\subset \Om_2$ constructed in the proof of Lemma \ref{h-body}). Therefore when $x\in l_{\ga}$ radially tends to a fixed point $p$ (in $\fix (\ga)\in S^2$) of such element $\ga$, its image $F(x)$ goes along that closed loop $F([x, \ga(x)])\subset\Om_1$ because $F(\ga(x))=\rho(\ga)(F(x))=F(x)$. Immediately it implies that the radial limit $\lim_{x\to p} F(x)$ does not exist. This shows that fixed points of any element $\ga\in F_3\subset \Ga$ (or its conjugates) are essential (topological) singularities of our quasiregular mapping $F$. So our quasiregular mapping $F$ has no continuous extension to the subset $S_*\subset S^2$ which is a dense subset of the unit sphere $S^2=\partial B^3\subset S^3$.
\end{cnstr}


  Dept of Math., University of Oklahoma, Norman, OK 73019, USA

   e-mail: apanasov\char`\@ ou.edu

\end{document}